\documentstyle[11pt,amssymb,amstex]{amsart}
\numberwithin{equation}{section} \textwidth 140mm \textheight 220mm
\def\cb{\mathcal{B}}

\def\cm{\mathcal{ M}}

\def\bn{{\mathbb N}}

\def\br{{\mathbb R}}

\def\l{\lambda} 

\def\m{\mu}
\def\p{\psi}

\def\s{\sigma} 
\def\t{\tau}
\def\f{\varphi}

\def\w{\omega}

\def\id{{\bf 1}\!\!{\rm I}}

\def\xb{{\mathbf{x}}}

\def\h{\frak{h}}

\newtheorem{thm}{Theorem}[section]

\newtheorem{cor}[thm]{Corollary}
\newtheorem{prop}[thm]{Proposition}
\newtheorem{defin}[thm]{Definition}
\newtheorem{rem}[thm]{Remark}

\def\o{\otimes}
\begin{document}

\title{ON MARGINAL MARKOV PROCESSES OF QUANTUM QUADRATIC STOCHASTIC PROCESSES}

\author{Farrukh MUKHAMEDOV}

\address{
Department of Computational \& Theoretical Sciences\\
Faculty of Science, International Islamic University Malaysia\\
P.O. Box, 141, 25710, Kuantan\\
Pahang, Malaysia $^*$ \\E-mail: far75m@@yandex.ru; farrukh\_m@@iiu.edu.my\\
}

\begin{abstract}
In the paper it is defined two marginal Markov processes on von
Neumann algebras $\cm$ and $\cm\o\cm$, respectively, corresponding
to given quantum quadratic stochastic process (q.q.s.p.). It is
proved that such marginal processes uniquely determines the q.q.s.p.
Moreover, certain ergodic relations between them are established as
well.
\end{abstract}

\keywords{quantum quadratic stochastic process; marginal Markov
process; ergodic principle}

\maketitle

\section{Introduction}

It is known that Markov processes are well-developed field of
mathematics which have various applications in physics, biology and
so on. But there are some physical models which cannot be described
by such processes. One of such models is a model related to
population genetics. Namely,  this model is described by quadratic
stochastic processes (see \cite{L,G1,SG1}). To define it, we denote
\begin{eqnarray*}
& & \ell^{1}=  \{x=(x_{n}):\|
x\|_{1}=\overset{\infty}{\underset{n=1}{\sum}}|x_{n}|<\infty ;\
x_{n}\in \mathbb{R}\},\\
&& S^{\infty} =\{x\in \ell^{1}:x_{n}\geq  0 ;\ \|x\|_{1}=1\}.
\end{eqnarray*}
Hence this process is defined as follows (see \cite{G1},\cite{SG1}):
Consider a family of functions $\{ p_{ij,k}^{[s,t]}: i,j,k\in\bn,\
s,\ t\in \br^{+},\ t-s \geq 1\}$. Such a family is said to be {\it
quadratic stochastic process (q.s.p.)} if for fixed $s,t\in\br_+$
it satisfies the following conditions:\\
\begin{enumerate}
   \item[(i)] $p_{ij,k}^{[s,t]}=p_{ji,k}^{[s,t]}$ for any $i,j,k\in\mathbb{N}$.
   \item[(ii)] $p_{ij,k}^{[s,t]}\geq 0$ and
$\overset{\infty}{\underset{k=1}{\sum }}p_{ij,k}^{[s,t]}=1$ for any
$i,j,k\in \mathbb{N}$.
   \item[(iii)] An analogue of Kolmogorov-Chapman equation; here there are two
variants: for the initial point $x^{(0)}\in S$,
$x^{(0)}=(x^{(0)}_{1},x^{(0)}_{2},\cdots)$ and $s<r<t$ such that $t-r\geq 1$,$r-s\geq 1$\\
(iii$_{A}$)
$$
p_{ij,k}^{[s,t]}=\overset{\infty}{\underset{m,l=1}{\sum
}}p_{ij,m}^{[s,r]}p_{ml,k}^{[r,t]}x^{(r)}_{l}, $$ where
$x^{(r)}_{k}$ is given by
$$
x^{(r)}_{k}=\overset{\infty}{\underset{i,j=1}{\sum
}}p_{ij,k}^{[0,r]}x^{(0)}_{i}x^{(0)}_{j};
$$
(iii$_{B}$)
$$
p_{ij,k}^{[s,t]}=\overset{\infty}{\underset{m,l,g,h=1}{\sum
}}p_{im,l}^{[s,r]}p_{jg,h}^{[s,r]}p_{lh,
k}^{[r,t]}x^{(s)}_{m}x^{(s)}_{g}.$$
\end{enumerate}

It is said that the q.s.p. $\{p_{ij,k}^{s,t}\}$ is of {\it type (A)
or (B)} if it satisfies the  fundamental equations (iii$_{A}$) or
(iii$_{B}$), respectively. In this definition the functions
$p_{ij,k}^{[s,t]}$ denotes the probability that under the
interaction of the elements $i$ and $j$ at time $s$ the element $k$
comes into effect at time $t$. Since for physical, chemical and
biological phenomena a certain time is necessary for the realization
of an interaction, we shall take the greatest such time to be equal
to 1 (see the Boltzmann model \cite{J} or the biological model
\cite{L}). Thus the probability $p_{ij,k}^{[s,t]}$ is defined for
$t-s\geq 1$.

It should be noted that such quadratic stochastic processes are
related to the notion of a quadratic stochastic operator, which was
introduced in \cite{B}, in the same way as Markov processes are
related to linear transformations (i.e. Markov operators). The
problem of studying the behaviour of trajectories of quadratic
stochastic operators was stated in \cite{U}. The limit behaviour and
ergodic properties of trajectories of such operators were studied in
( see for example \cite{Ke,L,Mak,V}.

 We note that quadratic stochastic processes describe physical systems defined above, but
they do not occupate the cases in quantum level. So, it is naturally
to define a concept of quantum quadratic processes. In \cite{GM,M}
quantum (noncommutative) quadratic stochastic processes (q.q.s.p.)
were defined on a von Neumann algebra and studied certain ergodic
properties ones. In \cite{GM} it is obtained necessary and
sufficient conditions for the validity of the ergodic principle for
q.q.s.p. From the physical point of view this means that for
sufficiently large values of time the system described by such a
process does not depend on the initial state of the system. It has
been found relations between quantum quadratic stochastic processes
and non-commutative Markov processes. In \cite{M} an expansion of
q.q.s.p. into a so-called fibrewise Markov process is given, and it
is proved that such an expansion uniquely determines the q.q.s.p. As
an application, it is given a criterion (in terms of this expansion)
for the q. q. s. p. to satisfy the ergodic principle. Using such a
result, it is proved that a q.q.s.p. satisfies the ergodic principle
if and only if the associated Markov process satisfies this
principle. It is natural to ask: is the defined Markov process
determines the given q.q.s.p. uniquely, or how many Markov processes
are needed to uniquely determine the q.q.s.p.? In this paper we are
going to affirmatively solve the problem. Namely, we shall show that
there two non-stationary Markov processes defined on different von
Neumann algebras $\cm$ and $\cm\o\cm$, respectively, called marginal
Markov processes, which uniquely determine the given to quantum
quadratic stochastic process. Such a description allows us to
investigate other properties of q.q.s.p. by means of Markov
processes. Moreover, certain ergodic relations between them are
established as well.

\section{Preliminaries}

Let $\cm$ be a von Neumann algebra acting on a Hilbert space $H$.
The set of all continuous (resp. ultra-weak continuous) functionals
on $\cm$ is  denoted by $\cm^*$ (resp. $\cm_{*}$), and put
$\cm_{*,+} = \cm_*\cap \cm^*_+$, here $\cm^*_+$ denotes the set of
all positive linear functionals. By $\cm\o\cm$ we denote tensor
product of $\cm$ in into itself. By $S$ and $S^2$ we denote the set
of all normal states on $\cm$ and $\cm\o \cm$ respectively. Recall a
mapping $U:\cm\o\cm\to\cm\o\cm$ is a linear operator such that
$U(x\o y)=y\o x$ for all $x,y\in \cm$. Let $\f\in S$ be a fixed
state. We define the conditional expectation operator
$E_{\f}:{\cm}\o{\cm}\to {\cm}$ on elements $a\o b$, $a,b\in{\cm}$ by
\begin{equation}\label{expect}
E_{\f}(a\o b)=\f(a)b
\end{equation} and extend it by linearity and continuity to ${\cm}\o{\cm}$.
Clearly, such an operator is completely positive and
$E_{\f}\id_{{\cm}\o{\cm}}=\id_{\cm}$ (more details on von Neumann
algebras we refer to \cite{T}).

Consider a family of linear operators
$\{P^{s,t}:{\cm}\to{\cm}\o{\cm}, s,t\in \br_+, t-s\geq 1 \}$.
\begin{defin}\label{qqsp} We say that a pair $(\{P^{s,t}\},\w_0)$, where $\w_0\in S$ is an initial state,
forms a {\it quantum quadratic stochastic process (q.q.s.p)}, if
every operator $P^{s,t}$ is ultra-weakly continuous and  the
following conditions hold:
\begin{enumerate}
\item[(i)] each operator $P^{s,t}$ is a unital completely positive mapping with
$UP^{s,t}=P^{s,t}$;

\item[(ii)] An analogue of Kolmogorov-Chapman equation is
satisfied: for initial state $\w_{0}\in S$ and arbitrary numbers
$s,\t,t\in \br_+$ such that $\t-s\geq 1, \ t-\t\geq  1$ one has

\item[(ii)$_{A}$]  $P^{s,t}x=P^{s,\t}(E_{\w_{\t}}(P^{\t,t}x)), \ \
x\in {\cm}$

\item[(ii)$_{B}$] $P^{s,t}x=E_{\w_{s}}P^{s,\t}\o
E_{\w_{s}}P^{s,\t}(P^{\t,t}x), \ \ x\in \cm$,
\end{enumerate}
where $\w_{\t}(x)=\w_{0}\o \w_{0}(P^{0,\t}x), \ x\in \cm$.
\end{defin}

If for q.q.s.p. the  fundamental equations (ii)$_A$ or (ii)$_B$ are
held  then we say that q.q.s.p.  has {\it type (A)} or {\it type
(B)}, respectively.

\begin{rem} By using the q.q.s.p., we can specify a law of interaction of states. For
$\f,\p\in S$, we set
$$
V^{s,t}(\f,\p)(x)=\f\o\p(P^{s,t}x), \ \qquad x\in\cm.
$$
This equality gives a rule according to which the state
$V^{s,t}(\f,\p)$ appears at time $t$ as a result of the interaction
of states $\f$ and $\p$ at time $s$. From the physical point of
view, the interaction of states can be explained as follows:
Consider two physical systems separated by a barrier and assume that
one of these systems is in the state $\f$. and the other one is in
the state $\p$. Upon the removal of the barrier, the new physical
system is in the state $\f\o\p$ and, as a result of the action of
the operator $P^{s,t}$, a new state is formed. This state is exactly
the result of the interaction of the states $\f$ and $\p$.
\end{rem}

\begin{rem} If $\cm$ is an $\ell^{\infty}$, i.e., $\cm = \ell^\infty$,
then a q.q.s.p. $(\{P^{s,t}\},\w_0)$ defined on $\ell^\infty$
coincides with a quadratic stochastic process. Indeed, we set
$$
p^{[s,t]}_{ij,k}=P^{s,t}(\chi_{\{k\}})(i,j), \ \ i,j,k\in\bn,
$$
where $\chi_A$ is the indicator of the set $A\in\cb$. Then, by
Definition \ref{qqsp}, the family of functions $p^{[s,t]}_{ij,k}$
forms a quadratic stochastic process.

Conversely, if we have a quadratic stochastic process
$(\{p^{[s,t]}_{ij,k}\},\m^{(0)})$ then we can define a quantum
quadratic stochastic process on $\ell^\infty$ as follows:
$$
(P^{s,t}{\mathbf{f}})(i,j)=\sum_{k=1}^\infty f_kp^{[s,t]}_{ij,k}, \
\ \mathbf{f}=\{f_k\}\in\ell^\infty.
$$
As the initial state, we take the following state
$$
\w_0({\mathbf{f}})=\sum_{k=1}^\infty f_k\m^{(0)}_k.
$$
One can easily check the conditions of Definition \ref{qqsp} are
satisfied. Thus, the notion of quantum quadratic stochastic process
generalizes the notion of quadratic stochastic process.
\end{rem}

\begin{rem} Certain examples of q.q.s.p were given in \cite{GM}.
\end{rem}

Let $(\{P^{s,t}\},\w_0)$ be a q.q.s.p. Then by $P^{s,t}_{*}$ we
denote the linear operator mapping from $(\cm\o \cm)_{*}$ into
$\cm_{*}$ given by
\begin{equation*}\label{7cqqsp}
P^{s,t}_{*}(\f)(x)=\f(P^{s,t}x), \ \ \f\in ({\cm}\o {\cm})_{*}, \ \
{x}\in {\cm}.
\end{equation*}

\begin{defin}\label{7ep} A q.q.s.p. $(\{P^{s,t}\},\w_0)$ is said to satisfy
{\it the ergodic  principle}, if  for every $\f,\p\in S^2$ and
$s\in\br_+$
$$
\lim_{t\to\infty}\|P^{s,t}_{*}\f-P^{s,t}_{*}\p\|_{1}=0,
$$
where $\|\cdot\|_1$ is the norm on $\cm^*$.
\end{defin}

Let us note that Kolmogorov was the first who introduced the concept
of an ergodic principle for Markov processes (see, for example,
\cite{K}). For quadratic stochastic processes such a concept was
introduced and studied in \cite{SG2,GHM}.

\section{Marginal Markov Processes and Ergodic Principle}

In this section we are going to consider relation between q.q.s.p.
and Markov processes.

First recall that a family $\{Q^{s,t}:\cm\to\cm ,  s,t\in\br_+,\
t-s\geq 1\}$ of unital completely positive operators is called {\it
Markov process} if
\begin{equation*}
Q^{s,t}=Q^{s,\t}Q^{\t,t}
\end{equation*}
holds for any $s,\t,t\in\br_+$ such that $t-\t\geq 1$, $\t-s\geq 1$.

A Markov process  $\{Q^{s,t}\}$ is said to satisfy the {\it ergodic
principle} if for every $\f,\p\in S$ and $s\in\br_+$ one has
$$
\lim_{t\to\infty}\|Q^{s,t}_*\f-Q^{s,t}_*\p\|_{1}=0.
$$

Let $(\{P^{s,t}\},\w_0)$ be a q.q.s.p. Then define a new process
$Q_P^{s,t}:{\cm}\to {\cm}$ by
\begin{equation}\label{7Markov-q2}
Q_P^{s,t}=E_{\w_s}P^{s,t}.
\end{equation}

Then according to Proposition 4.3 \cite{M} $\{Q^{s,t}_P\}$ is a
Markov process associated with q.q.s.p. It is evident that the
defined process satisfies the ergodic principle if the q.q.s.p.
satisfies one. An interesting question is about the converse.
Corollary 4.4 \cite{M}  states the following important result:

\begin{thm}\label{7ep-mp1}  Let $(\{P^{s,t}\},\w_0)$ be a q.q.s.p. on a von Neumann algebra $\cm$ and let
$\{Q^{s,t}_P\}$ be the corresponding Markov process.  Then the
following conditions are equivalent
\begin{enumerate}
\item[(i)] $(\{P^{s,t}\},\w_0)$ satisfies the ergodic principle;
\item[(ii)] $\{Q^{s,t}_P\}$ satisfies the ergodic principle;
 \item[(iii)] There is a number $\l\in[0,1)$ such that, given any
 states $\f,\p\in S^2$ and a number $s\in\br_+$ one has
\begin{equation*}
\|Q^{s,t}_{P,*}\f-Q^{s,t}_{P,*}\p\|_1\leq \l\|\f-\p\|_1
\end{equation*}
for at least one $t\in\br_+$.
\end{enumerate}
\end{thm}

\subsection{Case type A}

In this subsection we assume that q.q.s.p. $(\{P^{s,t}\},\w_0)$ has
type (A).

Now define another process $\{H^{s,t}_P:\cm\o \cm\to\cm\o \cm, \
s,t\in\br_+,\ t-s\geq 1\}$ by
\begin{eqnarray}\label{7H}
H^{s,t}_P\xb=P^{s,t}E_{\w_t}\xb,\ \qquad  \ \xb\in \cm\o \cm.
\end{eqnarray}

It is clear that every $H^{s,t}_P$ is a unital completely positive
operator.  It turns out that  $\{H^{s,t}\}$ is Markov process.
Indeed, using (ii)$_A$ of Def. \ref{qqsp} one has
\begin{eqnarray*}
H^{s,t}_P\xb&=&P^{s,t}E_{\w_t}\xb=P^{s,\t}E_{\w_\t}(P^{\t,t}E_{\w_t}\xb)=H^{s,\t}_PH^{\t,t}_P\xb,
\end{eqnarray*}
which is the assertion.

The defined two Markov processes $Q^{s,t}_P$ and $H^{s,t}_P$ are
related with each other by the following equality
$$
E_{\w_s}(H^{s,t}_P\xb)=E_{\w_s}(P^{s,t}(E_{\w_t}(\xb))=Q^{s,t}_P(E_{\w_t}(\xb))
$$
for every $\xb\in \cm\o \cm$. Moreover, $H^{s,t}_P$ has the
following properties
\begin{eqnarray}\label{7HE}
H^{s,t}_P\xb=P^{s,t}(E_{\w_t}(\xb))=P^{s,t}E_{\w_t}E_{\w_t}(\xb)=H^{s,t}_P(E_{\w_t}(\xb)\o\id)
\end{eqnarray}
\begin{eqnarray*}
UH^{s,t}_P=H^{s,t}_P, \ \ H^{s,t}_P(x\o\id)=P^{s,t}x, \ \  x\in \cm.
\end{eqnarray*}
from \eqref{7HE} one gets $ H^{s,t}_P(\id\o x)=\w_t(x)\id\o\id$.
Here we can represent
\begin{eqnarray*}
\w_t(x)=\w_0\o\w_0(P^{0,t}x)=\w_0(Q^{0,t}_Px),\\
\w_t(x)=\w_0\o\w_0(H^{0,t}_P(x\o\id)).
\end{eqnarray*}

Now we are interested in the following question: can such kind of
two Markov processes (i.e. with above properties) determine uniquely
a q.q.s.p.? To answer to this question we need to introduce some
notations.

Let  $\{Q^{s,t}:\cm\to\cm\}$ and $\{H^{s,t}:\cm\o \cm\to \cm\o
\cm\}$ be two Markov processes with an initial state $\w_0\in S$.
Denote
$$
\f_t(x)=\w_0(Q^{0,t}x), \ \ \ \p_t(x)=\w_0\o\w_0(H^{0,t}(x\o\id)).
$$
Assume that the given processes satisfy the following conditions:
\begin{enumerate}
\item[(i)] $UH^{s,t}=H^{s,t}$;\\
\item[(ii)] $E_{\p_s}H^{s,t}\xb=Q^{s,t}E_{\f_t}\xb$ for all $\xb\in\cm\o \cm$;\\
\item[(iii)] $ H^{s,t}\xb=H^{s,t}(E_{\p_t}(\xb)\o\id)$.
\end{enumerate}

First note that if we take $\xb=\id\o x$ in (iii) then we get
\begin{eqnarray}\label{7H1}
H^{s,t}(\id\o x)&=&H^{s,t}(E_{\p_t}(\id\o x)\o\id)\nonumber\\
&=&H^{s,t}(\p_t(x)\id\o\id)\nonumber\\
&=&\p_t(x)\id\o\id
\end{eqnarray}
Now from (ii) and \eqref{7H1} we have
\begin{eqnarray}\label{7H2}
E_{\p_s}H^{s,t}(\id\o x)&=&E_{\p_s}(\p_t(x)\id\o\id)\nonumber
\\
&=&\p_t(x)\id\nonumber\\
&=&Q^{s,t}E_{\f_t}(\id\o x)\nonumber\\
&=&\f_t(x)\id.
\end{eqnarray}
This means that $\f_t=\p_t$, therefore in the sequel we denote
$\w_t:=\f_t=\p_t.$

Now we are ready to formulate the result.

\begin{thm}\label{7HPQ} Let  $\{Q^{s,t}\}$ and $\{H^{s,t}\}$ be two Markov
Processes with (i)-(iii).  Then by the equality
$P^{s,t}x=H^{s,t}(x\o\id)$ one defines a q.q.s.p. of type (A).
Moreover, one has
\begin{enumerate}
\item[(a)] $P^{s,t}=H^{s,\t}P^{\t,t}$ for any $\t-s\geq 1$,
$t-\t\geq 1$,\\
\item[(b)] $Q^{s,t}=E_{\w_s}P^{s,t}$.
\end{enumerate}
\end{thm}

\begin{pf} We have to check only the condition (ii)$_A$ of Def. \ref{qqsp}. Take any
$\t-s\geq 1$, $t-\t\geq 1$. Then using the assumption (iii) we have
\begin{eqnarray*}
P^{s,\t}E_{\w_\t}(P^{\t,t}x)&=&H^{s,\t}(E_{\w_\t}H^{\t,t}(x\o\id)\o\id)\\
&=&H^{s,\t}H^{\t,t}(x\o\id)\\
&=&H^{s,t}(x\o\id)\\
&=&P^{s,t}x, \ \ x\in M.
\end{eqnarray*}
From Markov property of $H^{s,t}$ we immediately get (a).

If we put $\xb=x\o\id$ to (iii) then from \eqref{expect} one finds
$$
E_{\w_s}P^{s,t}x=E_{\w_s}H^{s,t}(x\o\id)=Q^{s,t}E_{\w_t}(x\o\id)=Q^{s,t}x.
$$
This completes the proof.
\end{pf}

These two $\{Q^{s,t}\}$ and $\{H^{s,t}\}$ Markov processes are
called {\it marginal Markov processes} associated with q.q.s.p.
$\{P^{s,t}\}$. So, according to Theorem \ref{7HPQ} the marginal
Markov processes uniquely define q.q.s.p.

Now define an other process $\{Z^{s,t}:\cm\o \cm\to \cm\o \cm\}$ by
\begin{eqnarray}\label{7Z}
Z^{s,t}\xb=&E_{\w_s}H^{s,t}(\xb)\o\id, \ \ \xb\in \cm\o \cm.
\end{eqnarray}
From (ii) one gets $Z^{s,t}\xb=Q^{s,t}E_{\w_t}\xb\o\id$. In
particular,
\begin{eqnarray*}\label{7Z1}
 && Z^{s,t}(x\o\id)=Q^{s,t}x\o\id,\\
&& Z^{s,t}(\id\o x)=\w_t(x)\id\o\id.
\end{eqnarray*}

\begin{prop}\label{7Z2} The process $\{Z^{s,t}\}$ is a Markov one.
\end{prop}

\begin{pf} Take any $\t-s\geq 1$, $t-\t\geq 1$. Then using the assumption (iii) and markovianity of $H^{s,t}$
we have
\begin{eqnarray*}
Z^{s,\t}Z^{\t,t}\xb&=&E_{\w_s}H^{s,\t}(E_{\w_\t}H^{\t,t}(\xb)\o\id)\o\id\\
&=&E_{\w_s}H^{s,\t}H^{\t,t}(\xb)\o\id\\
&=&E_{\w_s}H^{s,t}(\xb)\o\id\\
&=&Z^{s,t}\xb,
\end{eqnarray*}
for every $\xb\in \cm\o \cm$, which is the assertion.
\end{pf}

\begin{rem}\label{7ZH} Consider a q.q.s.p. $(\{P^{s,t}\},\w_0)$ of type (A). Let
$H^{s,t}$, $Z^{s,t}$ be the associated Markov processes.  Take any
$\f\in S^2$ then from \eqref{7H} with taking into account
\eqref{expect}, one concludes that
\begin{eqnarray}\label{7ZH1}
\f(H^{s,t}\xb)=P^{s,t}_*\f(E_{\w_t}(\xb))=P^{s,t}_*\f\o\w_t(\xb),
\end{eqnarray}
for any $\xb\in\cm\o\cm$.

Similarly, using \eqref{7Z}, for $Z^{s,t}$ we have
\begin{eqnarray}\label{7ZH2}
Z^{s,t}_*(\s\o\p)=\p(\id)P^{s,t}_*(\s\o\w_s)\o\w_t,
\end{eqnarray}
for every $\s,\p\in \cm_*$.
\end{rem}

From Theorem \ref{7ep-mp1} and using \eqref{7ZH1},\eqref{7ZH2} one
can prove the following

\begin{cor}\label{7ep-mp2}  Let $(\{P^{s,t}\},\w_0)$ be a q.q.s.p. of type (A) on $\cm$ and let
$\{Q^{s,t}\}$, $\{H^{s,t}\}$ be its marginal processes. Then the
following conditions are equivalent
\begin{enumerate}
\item[(i)] $(\{P^{s,t}\},\w_0)$ satisfies the ergodic principle;
\item[(ii)] $\{Q^{s,t}\}$ satisfies the ergodic principle;
\item[(iii)] $\{H^{s,t}\}$ satisfies the ergodic principle;
\item[(iii)] $\{Z^{s,t}\}$ satisfies the ergodic principle;
\end{enumerate}
\end{cor}

\subsection{Case type B}

Now suppose that a q.q.s.p. $(\{P^{s,t}\},\w_0)$ has type (B).

Like \eqref{7H} let us define a process $\h^{s,t}_P:\cm\o
\cm\to\cm\o \cm$ by
\begin{eqnarray}\label{7h}
{\h}^{s,t}_P\xb=P^{s,t}E_{\w_t}\xb,\ \qquad  \ \xb\in \cm\o \cm.
\end{eqnarray}

The defined process $\{{\h}^{s,t}_P\}$ is not Markov, but satisfies
other equation. Namely, using (ii)$_B$ of Def. \ref{qqsp} and
\eqref{7Markov-q2}  we get
\begin{eqnarray*}
{\h}^{s,t}_P\xb&=&E_{\w_s}P^{s,\t}\o
E_{\w_s}P^{s,\t}(P^{\t,t}E_{\w_t}\xb) =Q^{s,\t}_P\o
Q^{s,t}_P({\h}^{\t,t}_P\xb),
\end{eqnarray*}
where $\xb\in\cm\o\cm$.

Note that the process $\{{\h}^{s,t}_P\}$ has the same properties
like $\{H^{s,t}_P\}$.

Similarly to Theorem \ref{7HPQ} we can formulate the following

\begin{thm}\label{7hPQ} Let  $\{Q^{s,t}\}$ be a Markov process and $\{{\h}^{s,t}\}$ be another
processes, which satisfy (i)-(iii) and
\begin{eqnarray}\label{7hh}
{\h}^{s,t}=Q^{s,\t}\o Q^{s,\t}\circ \h^{\t,t}
\end{eqnarray}
for any $\t-s\geq 1$, $t-\t\geq 1$. Then by the equality
$P^{s,t}x={\h}^{s,t}(x\o\id)$ one defines a q.q.s.p. of type (B).
Moreover, one has $Q^{s,t}=E_{\w_s}P^{s,t}$.
\end{thm}

\begin{pf} We have to check only the condition (ii)$_B$. Note
that the assumption (iii) implies that
$$
E_{\w_s}{\h}^{s,t}(x\o\id)=Q^{s,t}E_{\w_t}(\cdot\o\id)=Q^{s,t}x, \ \
x\in\cm.
$$
Using this equality with \eqref{7hh} for any $\t-s\geq 1$, $t-\t\geq
1$ we have
\begin{eqnarray*}
E_{\w_s}P^{s,\t}\o E_{\w_s}P^{s,\t}(P^{\t,t}x)&=&
E_{\w_s}{\h}^{s,\t}(\cdot\o\id)\o
E_{\w_s}{\h}^{s,\t}(\cdot\o\id)({\h}^{\t,t}(x\o\id))\\
&=&Q^{s,\t}\o Q^{s,\t}(\h^{\t,t}(x\o\id))\\
&=&{\h}^{s,t}(x\o\id)\\
&=&P^{s,t}x
\end{eqnarray*}
for any $x\in\cm$.

This completes the proof.
\end{pf}

These two processes $\{Q^{s,t}\}$ and $\{{\h}^{s,t}\}$ we call {\it
marginal processes} associated with q.q.s.p. $\{P^{s,t}\}$.

Now define process $\{z^{s,t}:\cm\o \cm\to \cm\o \cm\}$ by
\begin{eqnarray}\label{7z}
z^{s,t}\xb=&E_{\w_s}{\h}^{s,t}(\xb)\o\id, \ \ \xb\in \cm\o \cm.
\end{eqnarray}
For this process \eqref{7Z1} also holds.

\begin{prop}\label{7z1} The process $z^{s,t}$ is a Markov one.
\end{prop}

\begin{pf} First from Theorem \ref{7hPQ} and Propsition 4.3 \cite{M} we
conclude that
\begin{equation}\label{7z11}
E_{\w_s}Q^{s,t}=E_{\w_t}.
\end{equation}

Let us take any $s,\t,t\in\br_+$ with $\t-s\geq 1$, $t-\t\geq 1$.
Then from \eqref{7z} with \eqref{7hh},\eqref{7z11} one gets
\begin{eqnarray}\label{7z12}
z^{s,t}\xb&=&E_{\w_s}(Q^{s,\t}\o Q^{s,\t}({\h}^{\t,t}(\xb))\o\id\nonumber\\
&=&Q^{s,\t}E_{\w_s}Q^{s,\t}({\h}^{\t,t}(\xb))\o\id\nonumber\\
&=&Q^{s,\t}E_{\w_\t}({\h}^{\t,t}(\xb))\o\id.
\end{eqnarray}

On the other hand, using conditions (ii),(iii) we obtain
\begin{eqnarray*}
z^{s,\t}z^{\t,t}\xb&=&E_{\w_s}h^{s,\t}(E_{\w_\t}{\h}^{\t,t}(\xb)\o\id)\o\id\\
&=&E_{\w_s}h^{s,\t}{\h}^{\t,t}(\xb)\o\id\\
&=&Q^{s,\t}E_{\w_\t}{\h}^{\t,t}(\xb)\o\id
\end{eqnarray*}
for every $\xb\in \cm\o \cm$. This relation with \eqref{7z12} proves
the assertion.
\end{pf}

\begin{cor}\label{7ep-mp3}  Let $(\{P^{s,t}\},\w_0)$ be a q.q.s.p. of type (B) on $\cm$ and let
$\{Q^{s,t}\}$, $\{\h^{s,t}\}$ be its marginal processes. Then the
following conditions are equivalent
\begin{enumerate}
\item[(i)] $(\{P^{s,t}\},\w_0)$ satisfies the ergodic principle;
\item[(ii)] $\{Q^{s,t}\}$ satisfies the ergodic principle;
\item[(iii)] $\{{\h}^{s,t}\}$ satisfies the ergodic principle;
\item[(iii)] $\{z^{s,t}\}$ satisfies the ergodic principle;
\end{enumerate}
\end{cor}

\section*{Acknowledgement} The author thanks the MOHE grant
FRGS0308-91 and the MOSTI grant 01-01-08-SF0079.


\end{document}